\documentclass[11pt]{amsart}
\usepackage{graphicx}
\usepackage{amssymb}
\begin{document}
\noindent
{\LARGE \bf A Note on Lattice Coverings}

\bigskip\noindent
Fei Xue and Chuanming Zong

\bigskip
\noindent
School of Mathematical Sciences, Peking University,

\noindent
Beijing 100871, P. R. China.

\noindent
cmzong@math.pku.edu.cn

\vspace{0.6cm}
\begin{minipage}[b]{13cm}
{\bf Abstract.} Whenever $n\ge 3$, there is a lattice covering $C+\Lambda$ of $E^n$ by a centrally symmetric convex body $C$ such that $C$ does not contain any parallelohedron $P$ that $P+\Lambda $ is a tiling of $E^n$.
\end{minipage}

\vspace{0.8cm}
\noindent
{\Large\bf 1. Introduction}

\bigskip
Let $K$ denote an $n$-dimensional convex body and let $C$ denote a centrally symmetric one centered at the origin of $E^n$. In particular, let $P$ denote an $n$-dimensional parallelohedron.  In other words, there is a suitable lattice $\Lambda$ such that $P+\Lambda $ is a tiling of $E^n$.

In 1885, E.S. Fedorov \cite{fedo85} discovered that, in $E^2$ a parallelohedron is either a parallelogram or a centrally symmetric hexagon (Figure 1); in $E^3$ a parallelohedron can be and only can be a parallelotope, a hexagonal prism, a rhombic dodecahedron, an elongated octahedron, or a truncated octahedron (Figure 2).

\begin{figure}[ht]
\centering
\includegraphics[height=3.8cm,width=9.4cm,angle=0]{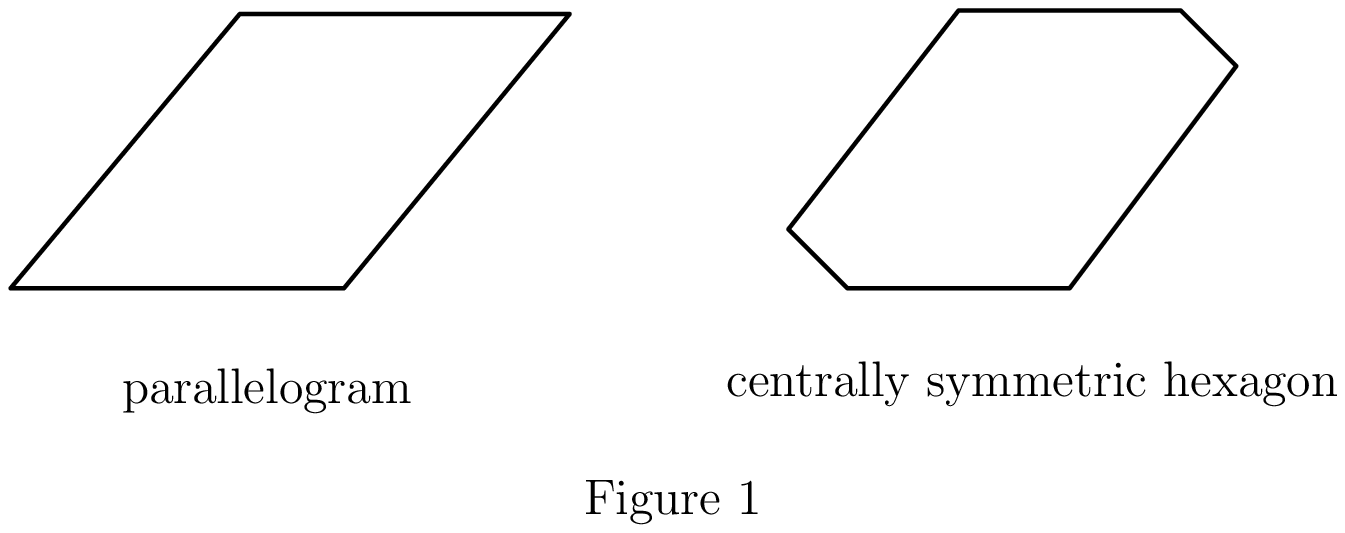}
\end{figure}

\begin{figure}[ht]
\centering
\includegraphics[height=10cm,width=11cm,angle=0]{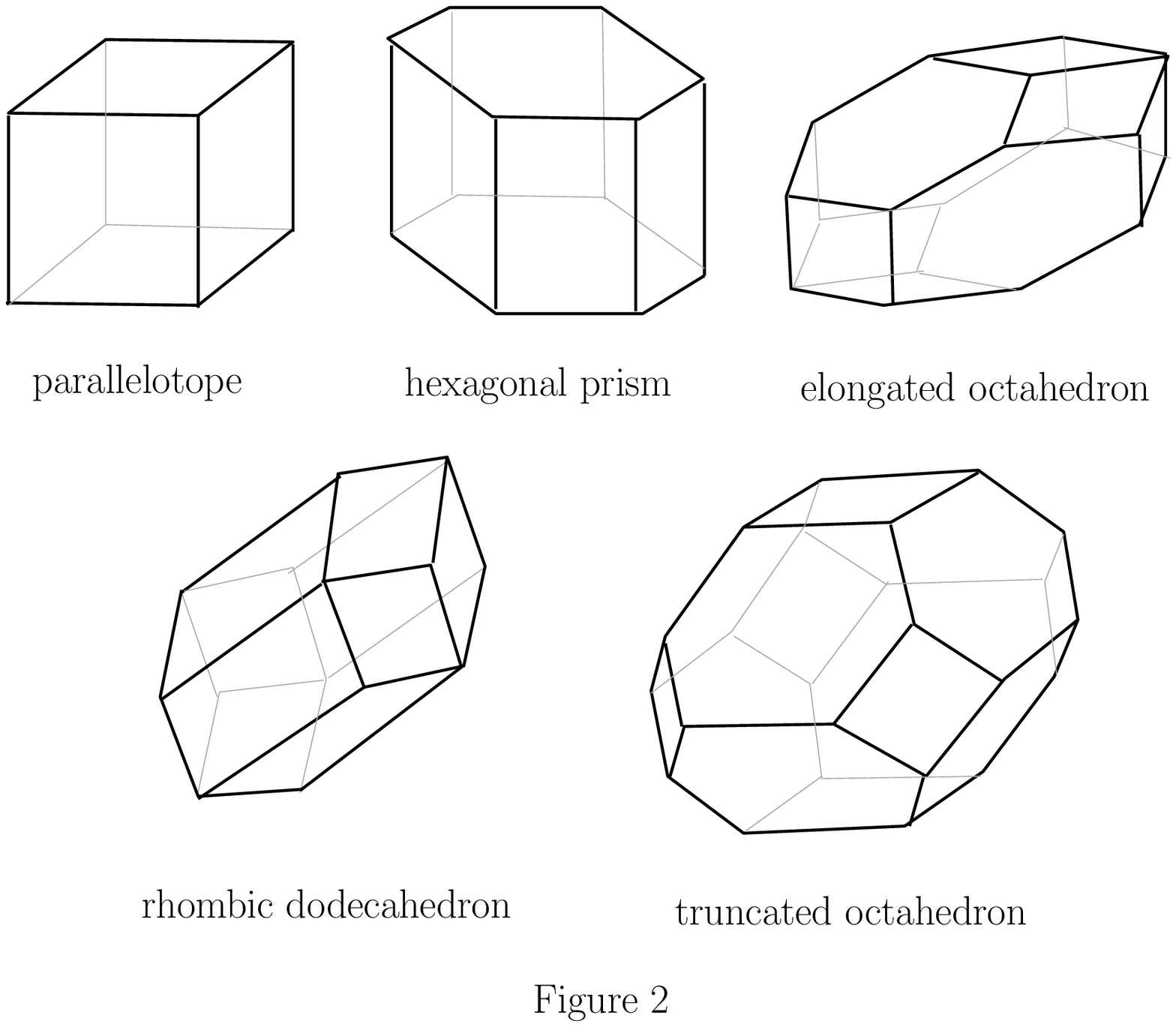}
\end{figure}

Let $\theta^t(K)$ denote the density of the thinnest translative covering of $E^n$ by $K$ and let $\theta^l(K)$ denote the density of the thinnest lattice covering of $E^n$ by $K$. For convenience,
let $B^n$ denote the $n$-dimensional unit ball and let $T^n$ denote the $n$-dimensional simplex with unit edges.
In 1939, Kerschner \cite{kers39} proved $\theta^t(B^2)= \theta^l(B^2)={{2\pi }/{\sqrt{27}}}.$
In 1946 and 1950, L. Fejes Toth \cite{feje46} and \cite{feje50} proved that $\theta^t(C)=\theta^l(C)\le {{2\pi}/ {\sqrt{27}}}$ holds for all two-dimensional centrally symmetric convex domains, where equality is attained precisely for the ellipses. In 1950, F\'ary \cite{fary50} proved that $\theta^l(K)\le {3/ 2}$ holds for all two-dimensional convex domains and the equality holds if and only if $K$ is a triangle. For more about coverings, we refer to \cite{bras05}, \cite{feje93} and \cite{pach95}.

If $K+\Lambda $ is a lattice covering of $E^2$, it can be easily shown that $K$ contains a centrally symmetric hexagon $H$ such that $H+\Lambda $ is a tiling of $E^2$. Therefore, let $\mathcal{H}$ denote the family of all centrally symmetric hexagons contained in $K$, we have
$$\theta^l(K)=\min_{H\in \mathcal{H}}{{{\rm vol} (K)}\over {{\rm vol}(H)}}.$$
This provides a practice method to determine the value of $\theta^l(K)$, in particular when $K$ is a polygon. Then, it is natural to raise the following problem in higher dimensions (see \cite{zong14}):

\medskip
\noindent
{\bf Problem 1.} {\it Whenever $K+\Lambda$ is a lattice covering of $E^n$, $n\ge 3$, is there always a parallelohedron $P$ satisfying both $P\subseteq K$ and $P+\Lambda $ is a tiling of $E^n$?}

\medskip
This note presents a counterexample to this problem.

\vspace{0.8cm}
\noindent
{\Large\bf 2. A Counterexample to Problem 1}

\bigskip
For convenience, we write $\alpha =\cos{\pi \over 3}$,  $\beta =\sin {\pi \over 3}$ and take $\gamma $ to be a small positive number. We note that $(1,0)$, $(\alpha , \beta )$, $(-\alpha , \beta )$, $(-1, 0)$, $(-\alpha , -\beta )$
and $(\alpha , -\beta )$ are the vertices of a regular hexagon. Let $C$ denote a three-dimensional centrally symmetric convex polytope as shown in Figure 3 with twelve vertices ${\bf v}_1=(1,0,1+\gamma )$, ${\bf v}_2=(\alpha , \beta , 1-\gamma )$, ${\bf v}_3=(-\alpha , \beta , 1+\gamma )$, ${\bf v}_4=(-1, 0 , 1-\gamma )$, ${\bf v}_5=(-\alpha , -\beta , 1+\gamma )$, ${\bf v}_6=(\alpha , -\beta , 1-\gamma )$, ${\bf v}_7=(1,0,-1+\gamma )$, ${\bf v}_8=(\alpha , \beta , -1-\gamma )$, ${\bf v}_9=(-\alpha , \beta , -1+\gamma )$, ${\bf v}_{10}=(-1, 0 , -1-\gamma )$, ${\bf v}_{11}=(-\alpha , -\beta , -1+\gamma )$ and ${\bf v}_{12}=(\alpha , -\beta , -1-\gamma )$, and let $\Lambda $ to be the lattice with a basis ${\bf a}_1=(1+\alpha , \beta , 0)$, ${\bf a}_2=(1+\alpha , -\beta , 0)$ and ${\bf a}_3=(0, 0, 2)$. In fact, $C$ can be obtained from an hexagonal prism of height $2(1+\gamma )$ by cutting off six tetrahedra, all of them are congruent to each others.

\begin{figure}[ht]
\centering
\includegraphics[height=5cm,width=4.5cm,angle=0]{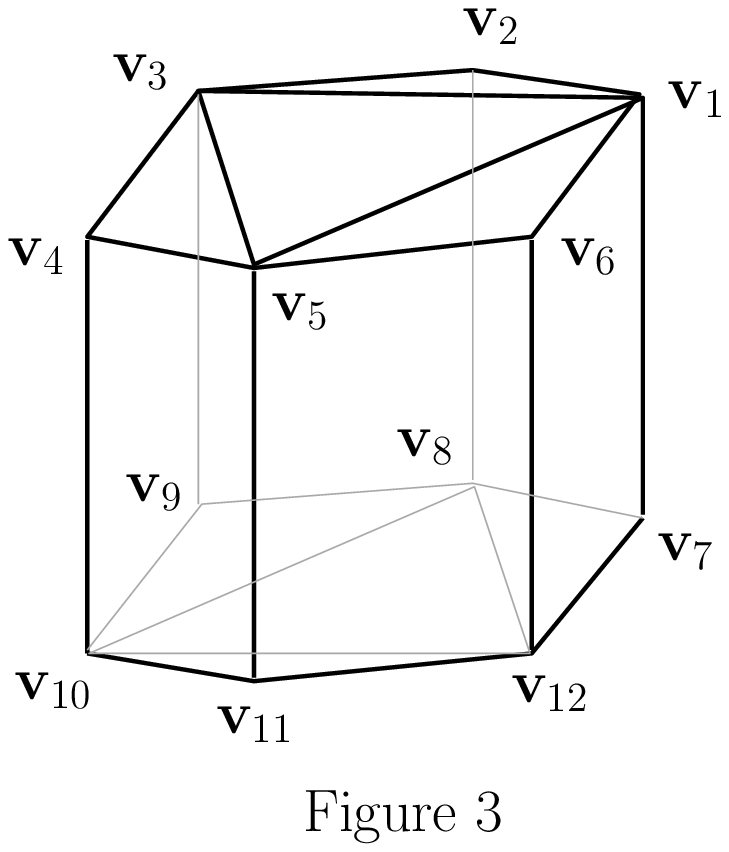}
\end{figure}

It can be easily verified that
$${\bf v}_i={\bf v}_{6+i}+{\bf a}_3$$
holds for all $i=1, 2, \ldots , 6$ and $C+\Lambda $ is a lattice covering of $E^n$. If $C$ contains a parallelohedron $P$ such that $P+\Lambda $ is a tiling of $E^3$, then $P$ must contain all the twelve vertices ${\bf v}_1$, ${\bf v}_2$, $\ldots $, ${\bf v}_{12}$ of $C$ and therefore $P=C$.
However, $C$ is apparently not a parallelohedron. Thus, $C+\Lambda $ is a counterexample to Problem 1 in
$E^3$.

\medskip
If $K$ is a counterexample to Problem 1 in $E^{n-1}$, defining $K'$ to be the cylinder over $K$,
one can easily show that $K'$ will be a counterexample to Problem 1 in $E^n$. Therefore, we have proved the following
result by explicit examples:

\medskip
\noindent
{\bf Theorem 1.} {\it Whenever $n\ge 3$, there is a lattice covering $C+\Lambda$ of $E^n$ by a centrally symmetric convex body $C$ such that $C$ does not contain any parallelohedron $P$ that $P+\Lambda $ is a tiling of $E^n$.}

\vspace{1cm}\noindent
{\bf Acknowledgements.}

\medskip
\noindent
This work is supported by 973 Programs 2013CB834201 and 2011CB302401, the National Science Foundation of China (No.11071003), and the Chang Jiang Scholars Program of China.

\vspace{0.6cm}
\bibliographystyle{amsplain}

\begin{thebibliography}{99}


\bibitem{bras05}P. Brass, W. Moser and J. Pach, {\it Research Problems in Discrete Geometry},
Springer-Verlag, New York, 2005.
\bibitem{fary50}I. F\'ary, Sur la densit\'e des r\'eseaux de domaines convexes, {\it Bull.
Soc. Math. France} {\bf 178} (1950), 152-161.
\bibitem{fedo85}E.S. Fedorov, Elements of the study of figures, {\it Zap. Mineral. Imper. S. Petersburgskogo Ob$\check{s}$$\check{c}$}, {\bf 21}(2) (1885), 1-279.
\bibitem{feje93}G. Fejes T\'oth and W. Kuperberg, Packing and
covering with convex sets, {\it Handbook of Convex Geometry} (P.M.
Gruber and J.M. Wills, eds.), North-Holland, Amsterdam 1993,
799--860.
\bibitem{feje46}L. Fejes Toth, Eine Bemerkung \"uber die Bedeckung der Eben durch Eibereiche mit
Mittelpunkt, {\it Acta Sci. Math. Szeged} {\bf 11} (1946), 93-95.
\bibitem{feje50}L. Fejes Toth, Some packing and covering theorems, {\it Acta Sci. Math. Szeged}
{\bf 12} (1950), 62-67.
\bibitem{kers39}R. Kerschner, The number of circles covering a set, {\it Amer. J. Math.}
{\bf 61} (1939), 665-671.
\bibitem{pach95}J. Pach and P.K. Agarwal, {\it Combinatorial Geometry}, John Wiley $\&$ Sons, 1995.
\bibitem{zong14}C. Zong, Minkowski bisectors, Minkowski cells, and lattice coverings, arXiv:1402.3395.

\end{thebibliography}

\end{document}